\documentclass[12pt]{article}

\usepackage{amsmath}
\usepackage{amssymb}
\usepackage{amsfonts}
\usepackage{lscape}

\begin{document}

\fontsize{12}{6mm}\selectfont
\setlength{\baselineskip}{2em}

$~$\\[.35in]
\newcommand{\dss}{\displaystyle}
\newcommand{\raro}{\rightarrow}
\newcommand{\be}{\begin{equation}}

\def\sech{\mbox{\rm sech}}
\def\sn{\mbox{\rm sn}}
\def\dn{\mbox{\rm dn}}
\thispagestyle{empty}

\begin{center}
{\Large\bf The Hodge-de Rham Decomposition Theorem }  \\    [2mm]
{\Large\bf And Some Applications Pertaining to Partial}  \\    [2mm]
{\Large\bf Differential Equation}  \\    [2mm]
\end{center}

\vspace{1cm}
\begin{center}
{\bf Paul Bracken}                        \\
{\bf Department of Mathematics,} \\
{\bf University of Texas,} \\
{\bf Edinburg, TX  }  \\
{78541-2999}
\end{center}

\vspace{3cm}
\begin{abstract}
The Hodge-de Rham Theorem is introduced and discussed.
This result has implications for the general study of
several partial differential equations.
Some propositions which which have applications to the proof
of this theorem are used to study some related 
results concerning a class of partial differential equation
in a novel way.
\end{abstract}

\vspace{2mm}
MSCs: 58A05, 58A10, 47J05.

\vspace{2mm}

\newpage

{\bf 1.} Let $M$ be a compact, orientable Riemannian
manifold. Let $E^{p} (M)$ denote the space of smooth
$p$-forms on $M$. An inner product can be defined on
the vector space $E^p (M)$ of $p$-forms on $M$ 
with all the usual properties {\bf [1,2]} by 
writing
$$
( \alpha, \beta ) = \int_M \, \alpha \wedge * \beta ,
\eqno(1)
$$
for $\alpha$, $\beta \in E^p (M)$ and $*$ is the 
Hodge operator {\bf [2]}. The corresponding
norm is denoted by $\left\| \alpha \right\|$. Suppose
that $\alpha$ and $\beta$ are forms of degree $p$
and $p+1$, respectively.
Stokes' Theorem can be used to show that
$$
\int_M \, d \alpha \wedge * \beta = 
(-1)^{p-1} \, \int_M \, \alpha \wedge d * \beta.
\eqno(2)
$$
Using (1), this can be put in the form
$$
( d \alpha, \beta ) = ( \alpha, \delta \beta ). 
\eqno(3)
$$
The operators $d$ and $\delta$ are adjoints or duals 
of each other {\bf [3,4]}.
Analogously, if $\beta$ is of degree $p-1$, then (2) gives
$( \alpha, d \beta ) = ( \delta \alpha, \beta )$. 
Thus in order that $\alpha$ be closed, it is necessary and
sufficient that it be orthogonal to all co-exact forms
of degree $p$. An operator called the Laplace-Beltrami 
operator is defined in terms of $d$ and $\delta$ as
$$
\Delta = d \delta + \delta d.
$$
The form $\omega \in E^p (M)$ is called harmonic if
$\Delta \omega =0$. Let $ H^p (M)$ denote the subspace
of $E^p (M)$ of harmonic forms on $M$.
From (3), it is clear that $\Delta$ is a self-adjoint
operator. It is clear that the harmonic forms of a
given degree will form a linear space.

Let $\beta$ be a $p$-form on a compact, orientable Riemannian
manifold $M$. If there is a $p$-form $\alpha$ such that
$\Delta \alpha = \beta$, then for a harmonic form
$\gamma \in H^p$
$$
( \beta, \gamma ) = ( \Delta \alpha, \gamma ) =
( \alpha , \Delta \gamma ) =0.
$$
Therefore, in order that a form $\alpha$ exist with
the property that $\Delta \alpha = \beta $, it is
necessary that $\beta$ be orthogonal to the subspace
$H^p$ of $E^p (M)$. This also turns out to be sufficient.
This reasoning might be thought to lead to the idea
that $E^p (M)$ can be partitioned into three distinct
subspaces which will be denoted $\Lambda^p_d$ and 
$\Lambda^p_{\delta}$ and $H^p$. These subspaces are 
respectively the ones which consist of those forms
which are exact, co-exact and harmonic, respectively.
These subspaces should be orthogonal in pairs, so forms
belonging to distinct subspaces are orthogonal.
This is in fact the case and when presented more
precisely is referred to as the Hodge-de Rham Decomposition Theorem.
It is the intention here to introduce this theorem and 
some associated results used in the proof. These results 
and techniques are then used to prove some related theorems
which are of great interest in a new way.

The Hodge-de Rham Theorem can be stated in a number 
of ways and two are now presented {\bf [1,5,6]}.

{\em HdR 1.} A regular differential form of degree $p$ may
be uniquely decomposed into a sum of the form
$$
\alpha = \alpha_d + \alpha_{\delta} + \alpha_{H},
$$
where $\alpha_d \in \Lambda^p_d $, $\alpha_{\delta} \in \Lambda^p_{\delta}$
and $\alpha_H \in H^p$.

{\em HdR 2.} For each integer $p$ with $0 \leq p \leq n$,
$H^p$ is finite dimensional and there is the direct sum
decomposition of the space $E^p (M)$ of smooth $p$-forms 
on $M$ given by
$$
E^p (M) = \Delta (E^p) \oplus H^p = 
d \delta ( E^p) \oplus \delta d (E^p) \oplus H^p.
$$
Consequently, the equation $\Delta \omega = \alpha$
has a solution in $E^p (M)$ if and only if the $p$-form
$\alpha$ is orthogonal to the space $H^p$.

{\bf 2.} The Hodge Theorem is connected to the
problem of finding necessary and sufficient conditions
for there to exist a solution $\omega$ for the equation
$\Delta \omega = \alpha$. Suppose that $\omega$ is a
solution of $\Delta \omega = \alpha$, then
$$
( \Delta \omega, \vartheta ) = ( \alpha, \vartheta ),
\eqno(4)
$$
for all $\vartheta \in E^p (M)$. Usually no
distinction is made between $\Delta$ and the  adjoint
$\Delta^{*}$, however, this distinction is important
in making the following definition.
Thus from (4), we can write
$$
( \omega, \Delta^{*} \vartheta ) = ( \alpha, \vartheta ).
\eqno(5)
$$
Equation (5) suggests that a solution of $\Delta \omega 
= \alpha$ can be regarded as a certain type of linear functional
on $E^p (M)$. In fact, $\omega$ determines a bounded linear
functional $L$ on $E^{p} (M)$ as follows
$$
L ( \beta ) = ( \omega, \beta ).
\eqno(5)
$$
Using (5), it follows that $L$ satisfies
$$
L ( \Delta^* \vartheta ) = ( \alpha, \vartheta ),
\eqno(7)
$$
for all $\vartheta \in E^p (M)$. Such a linear functional
is referred to as a {\em weak} solution of $\Delta \omega = \alpha$, 
or to state it another way, a weak solution of $\Delta \omega = \alpha$
is a bounded linear functional $L : E^p (M) \rightarrow \mathbb R$ 
such that (7) holds. Each ordinary solution determines 
a weak solution by (6). In fact, the converse is true, 
each weak solution determines an ordinary solution. The
Proposition which follows is fundamental for what will come.

{\bf Proposition 1.} Let $\alpha \in E^p (M)$ and 
suppose $L$ is a weak solution of  $\Delta \omega = \alpha$.
Then there exists $\omega \in E^p (M)$ such that
$$
L (\beta ) = ( \omega, \beta )
\eqno(8)
$$
for every $\beta \in E^p (M)$. Consequently $\Delta \omega = \alpha$.

This proposition as well as the next will be used in
the sequel.

{\bf Proposition 2.} Let $\{ \alpha_n \}$ be a sequence of 
smooth $p$-forms on $M$ such that $\left\| \alpha_n \right\| \leq c$ 
and $\left\| \Delta \alpha_n \right\| \leq c$ for all
$n$ and for some constant $c >0$. Then a subsequence of 
$\{ \alpha_n \}$ is a Cauchy sequence in $E^p (M)$.

Assuming the Propositions $1$ and $2$, a number of interesting results
can be developed including a proof of (HdR). It is the intention now 
to use these to carry this out.

{\bf Theorem 1.} For each integer $p$ with $0 \leq p \leq n$,
the space $H^p$ has finite dimension.

{\bf Proof:} If $H^p$ were not finite dimensional, then $H^p$ would
contain an infinite orthonormal sequence $\{ e_n \}$. Now
$\left\| e_n \right\| \leq 1$ and $ \left\| \Delta e_n \right\| =0$
since $e_n \in H^p$ for all $n$. By Proposition 2, a subsequence
of $\{ e_n \}$ is a Cauchy sequence in $E^p (M)$. Extract this sequence 
and call it simply $\{ e_n \}$. Let $\epsilon < 1$ be given, then
since this is Cauchy, there exists an $N ( \epsilon )$ such that
for $n,m > N ( \epsilon )$, $\left\| e_n - e_m \right\| < \epsilon$.
However, this is not possible since the $e_n$ are orthonormal so
$\left\| e_n - e_m \right\|^2 =2$.

{\bf Lemma 1.} Let $H$ denote the projection operator of
$E^p (M)$ onto $H^p$ so that $H ( \alpha )$ is the harmonic
part of $\alpha$, which will also be denoted $\alpha_H$.

$(i)$ For any $\alpha \in E^p (M)$, there exists a unique
$p$-form $H ( \alpha) \in H^p$ with the property that
$( \alpha, \gamma ) = ( H( \alpha), \gamma )$ for all
$\gamma \in H^p$. 
$(ii)$ For any $p$-form $\alpha \in E^p (M)$,
$H$ satisfies $ H (H (\alpha)) = H (\alpha)$.

{\bf Proof:} From Theorem 1, let $e_1, \cdots, e_m$ be
an orthonormal basis for $H^p$. Then an arbitrary form
$\alpha \in E^p (M)$ can be written uniquely in the form
$$
\alpha = \gamma + \sum_{i=1}^m \, ( \alpha, e_i ) \, e_i,
\eqno(9)
$$
where $\gamma \in (H^p)^{\perp}$, the subspace of $E^p (M)$
consisting of all elements orthogonal to $H^p$. Therefore,
$$
H (\alpha) = \sum_{i=1}^m \, ( \alpha, e_i) \, e_i =
\alpha - \gamma,
$$
and $( H(\alpha), \beta ) = ( \alpha - \gamma, \beta ) =
( \alpha, \beta )$ for all $\beta \in H^p$ since
$\gamma \in ( H^p )^{\perp} $. Moreover, since $H^p$
is a linear space $H ( H(\alpha)) = H (\alpha) -
H (\gamma ) = H ( \alpha )$, since $\gamma \in
(H^p)^{\perp}$.  $\clubsuit$

Consequently, Theorem 1 yields an orthogonal direct
sum decomposition of $E^p (M)$ which takes the form $E^p (M) =
(H^p)^{\perp} \oplus H^p$. The Hodge-de Rham 
Theorem then follows from this result by showing
that $ (H^p)^{\perp} = \Delta (E^p)$. Clearly,
$\Delta (E^p) \subset ( H^p)^{\perp}$ since if
$\omega \in E^p (M)$ and $\alpha \in H^p$, then
$$
( \Delta \omega, \alpha ) = ( \omega, \Delta \alpha ) =0.
$$
To complete the proof of the Hodge Theorem, the converse
of this must be shown, namely $H^{p} \subset \Delta (E^p)$.
To do so along the present lines, the following Lemma
is required, and it is stated without proof.

{\bf Lemma 2.} There exists a constant $k >0$ such that
$$
\left\| \gamma \right\| \leq k \left\| \Delta \gamma \right\|,
\qquad
\gamma \in (H^p )^{\perp}.
\eqno(10)
$$
Suppose $\alpha \in (H^p)^{\perp}$. Define a linear
functional $L$ on $\Delta (E^p)$ by setting
$$
L ( \Delta \gamma ) = ( \alpha, \gamma ),
\qquad
\gamma \in E^p (M).
$$
To show that $L$ is well defined, suppose that
$\Delta \gamma_1 = \Delta \gamma_2$, then
$\gamma_1 - \gamma_2 \in H^p$, which means that
$( \alpha, \gamma_1 - \gamma_2 ) =0$. Also $L$ is a
bounded linear functional on $\Delta (E^p)$, for
let $\gamma \in E^p (M)$ and $\psi = \gamma - H (\gamma)
\in (H^p)^{\perp}$. Using (10) in Lemma 2 and Cauchy-Schwarz
$$
| L (\Delta \gamma)| = | L ( \Delta \psi) |
= | ( \alpha, \psi)| \leq \left\| \alpha \right\| \left\|
\psi \right\| \leq k \left\| \alpha \right\|
\left\| \Delta \psi \right\|= k \left\| \alpha \right\|
\left\| \Delta \gamma \right\|.
$$
By the Hahn-Banach Theorem, $L$ extends to a
bounded linear functional on $E^p (M)$, and so
$L$ is a weak solution of $\Delta \omega = \alpha$.
By Proposition 1, it may be concluded that there
exists $\omega \in E^p (M)$ such that $\Delta \omega =
\alpha$. Therefore, it may be concluded that
$(H^p)^{\perp} = \Delta (E^p)$, and this gives $(HdR 2)$.

{\bf 3.} Consequently, for a given $p$-form $\alpha$,
there must exist a $p$-form $\beta$ which satisfies the
equation
$$
\Delta \beta = \alpha - H( \alpha).
\eqno(11)
$$
Thus from the Hodge Theorem and Lemma 1, it follows that
$\alpha_d + \alpha_{\delta} = \alpha - H(\alpha)$ and
$\beta$ can be taken such that $\Delta \beta = \alpha_d
+ \alpha_{\delta}$. Let $\beta_1$ and $\beta_2$ be any
two forms which satisfy (11). Subtracting the two
equations obtained from (11) which contain 
$\beta_1$ and $\beta_2$, there remains the equation
$\Delta ( \beta_1 - \beta_2 ) =0$. This result implies that
$\beta_1 - \beta_2 \in H^p$, and so modulo harmonic
forms, there exists a unique solution to (11). This solution
will be denoted here as $\beta = G (\alpha)$.
From Lemma 1, it is clear that $\beta$ can be
projected down onto $\Delta (E^p)$, and it is only this
part that need be substituted for $\beta$ in (11). Thus it suffices
to suppose that $G (\alpha)$ has no harmonic contribution.
Substituting $\beta$ into (11), we obtain
$$
\alpha = \Delta G (\alpha) + H ( \alpha),
\qquad
( G (\alpha), \sigma ) =0,
\eqno(12)
$$
for all $\sigma \in H^p$. This can now be summarized as Theorem 2.

{\bf Theorem 2.} Operator $G: E^p (M) \rightarrow (H^p)^{\perp}$
is defined by setting $G (\alpha)$ to be the unique solution of
(11) in $(H^p)^{\perp}$. For a given $p$-form $\alpha$, there
exists a $p$-form $G (\alpha)$ which satisfies differential
equation (11), namely $\Delta G (\alpha) = \alpha - H (\alpha)$,
and satisfies the associated condition given in (12).

It is important to realize that $G$ has important commutation 
properties which are revealed in the following Theorem.

{\bf Theorem 3.} $G$ commutes with any linear operator
which commutes with the Laplacian $\Delta$.

{\bf Theorem 4.} $(i)$ G is a self-adjoint operator
$$
( G \alpha, \beta ) = ( \alpha, G \beta ),
\eqno(13)
$$
for forms $\alpha$ and $\beta$ of the same degree $p$.

$(ii)$ $G$ is a Hermitean positive operator,
$$
( G \alpha, \alpha ) \geq 0,
\eqno(14)
$$
and equality holds if and only if $\alpha$ is harmonic.

{\bf Proof:} $(i)$ Based on Theorem 2, let $\alpha$ and $\gamma$
be $p$-forms such that $\Delta G \alpha = \alpha - \alpha_{H}$
and $\Delta G \gamma = \gamma - \gamma_{H}$. From orthogonality,
it follows that $( \alpha_{H}, \gamma) = ( \alpha, \gamma_{H})$.
Consequently, by means of the relations,
$$
( \Delta G \alpha, \gamma) = ( \alpha- \alpha_{H}, \gamma)
= ( \alpha, \gamma) - ( \alpha_{H}, \gamma ),
\quad
( \alpha, \Delta G \gamma ) = ( \alpha, \gamma - \gamma_{H})
= ( \alpha, \gamma) - ( \alpha, \gamma_{H}),
$$
it can be concluded that
$$
( \Delta G \alpha, \gamma ) = ( \alpha, \Delta G \gamma ).
$$
Since $\Delta$ and $G$ commute and $G$ is self-adjoint,
we can write the two expressions
$$
( \Delta G \alpha, \gamma) = ( G \alpha, \Delta \gamma),
\qquad
( \alpha, \Delta G \gamma) = ( \alpha, G \Delta \gamma).
$$
From these, it can be concluded that
$$
( G \alpha, \Delta \gamma) = ( \alpha, G \Delta \gamma ).
$$
Set $\beta = \Delta \gamma$ in this result and the required follows
$( G \alpha, \beta) = ( \alpha, G \beta )$. It suffices to
show this much, since the range of $G$ does not intersect
the harmonic subspace. If $\beta$ has a harmonic term, it will not
make a contribution, since it projects out giving zero by
orthogonality of these subspaces.

$(ii)$ Since $G \alpha$ has no component in $H^p$,
it follows by applying Theorem 2 that
$$
( G \alpha, \Delta G \alpha + H( \alpha) ) = ( G \alpha, \Delta G \alpha)
= ( \delta G \alpha , \delta G \alpha )
+ ( d G \alpha, d G \alpha ) \geq 0,
$$
since $( \sigma, \sigma ) \geq 0$ for any $p$-form
$\sigma$.

Suppose that $\alpha \in H^p$, then 
$( G \alpha, \alpha) =0$ by orthogonality. Suppose that
$( G \alpha, \alpha ) =0$ with $\alpha$, $G \alpha$ not
zero forms. Then since $G \alpha \in (H^p)^{\perp}$,
$\alpha$ must lie in the orthogonal subspace, $\alpha
\in H^p$ as discussed in Lemma 1, so $\alpha$ is purely harmonic.

${\bf 4.}$ Consider the Laplace-Beltrami operator $\Delta$
which acts on $p$-forms for some fixed $p$. A real
number $\lambda$ for which there is a $p$-form $\omega$
which is not identically zero such that $\Delta \omega =
\lambda \omega$ is called an eigenvalue of $\Delta$ and
the corresponding $p$-form $\omega$ an eigenfunction of
$\Delta$ corresponding to $\lambda$. The eigenfunctions 
corresponding to a fixed $\lambda$ form a subspace of 
$E^p (M)$, namely the eigenspace of $\lambda$.
Operator $G$ can also be used to discuss some of the
properties of the eigenvalues and their forms. It
is useful to start by reviewing some basic properties
of these objects which will be assumed and used later.

The eigenvalues of $\Delta$ are non-negative. Suppose
that $\alpha \in E^p (M)$ such that $\Delta \alpha = \lambda \alpha$.
It follows that $( \Delta \alpha, \alpha) = \lambda
(\alpha, \alpha)$, which can be written in the
equivalent form $( \delta \alpha, \delta \alpha)
+ ( d \alpha, d \alpha ) = \lambda \left\| \alpha \right\|^2$.
With the exception of $\lambda$, all the terms in this 
equality are known to be positive, therefore $\lambda$ 
must be positive as well.

The eigenspaces of $\Delta$ are finite dimensional.
Suppose this is not the case, then there exists an
eigenspace which contains an infinite orthonormal
sequence or basis, which will be called $\{ e_n \}$
such that $\left\| e_n \right\| =1$ and $\left\| 
\Delta e_n \right\| \leq \lambda$ for all $n$.
Then as in Theorem 1, this orthonormal sequence
would contain a Cauchy subsequence, which is not possible.

The eigenvalues have no finite accumulation point. 
Suppose the eigenvalues have such a value, $\mu$.
Let $\{ e_n \}$ be a set of orthonormalized
eigenfunctions of $\Delta$ associated to the
eigenvalues $\{ \lambda_j \}$, which accumulate to
$\mu$, with one selected each respective eigenspace.
Then $\left\| e_j \right\| =1$ is bounded for all
$j$ and $\left\| \Delta e_j \right\| = \lambda_j$.
Since the $\lambda_j$ accumulate to the value $\mu$,
$\sup_{j} \, \left\| \Delta e_j \right\| \leq C$.
Then there is a subsequence of the $\{ e_j \}$ which
is Cauchy, and again this is not possible.

The eigenfunctions corresponding to distinct
eigenvalues are orthogonal. Let $\alpha$, $\beta \in E^p (M)$
be eigenfunctions with $\Delta \alpha = \lambda \alpha$
and $\Delta \beta = \mu \beta$. Using the first of
these, $( \Delta \alpha, \beta) = \lambda ( \alpha,
\beta )$. Since $\Delta$ is self-adjoint, this is
equivalent to $( \alpha, \Delta \beta ) = \lambda
( \alpha, \beta )$. However, using the second equation gives
$( \alpha, \Delta \beta ) = \mu ( \alpha, \beta )$.
Combining these, it must be that $\lambda (\alpha, \beta )
= \mu ( \alpha, \beta )$. If $\lambda \neq \mu$,
then it follows that $( \alpha, \beta ) =0$, which
means $\alpha$ and $\beta$ are orthogonal.

{\bf 5.} The operator $G$ can be of use in discussing the
existence  of eigenvalues and some of their
properties. Clearly, zero is an eigenvalue of
$\Delta$ is and only if there are nontrivial harmonic 
$p$-forms on $M$. The eigenspace agrees exactly with the space
$H^p$ of harmonic forms on $M$.

{\bf Theorem 5.} The eigenvalues of $G$ acting on $(H^p)^{\perp}$
are the reciprocals of the nonzero eigenvalues of
$\Delta$ acting on $(H^p)^{\perp}$.

{\bf Proof:} Consider the restriction $\Delta :
(H^p)^{\perp} \rightarrow (H^p)^{\perp}$, which
matches the domain of $G$ as well. Any $\omega 
\in (H^p)^{\perp}$ satisfies $\Delta G \omega = \omega$,
which implies $G \Delta \omega = \omega$, since
$G$ and $\Delta$ commute. Suppose $\omega$ is an
eigenfunction of $\delta$ so $\Delta \omega = \lambda
\omega$ for some $\lambda \in \mathbb  R^{+}$.
By the opening remark, we have $ \omega = G \Delta \omega
= \lambda G \omega$. Since $\lambda \neq 0$ for such
$\omega$, this implies that $G \omega = \lambda^{-1} \omega$.

The following Theorem ensures that this is not an empty
statement.

{\bf Theorem 6.} The Laplacian $\Delta$ has a positive
eigenvalue and in fact a whole sequence of eigenvalues
which diverge to $+ \infty$.

{\bf Proof:} The argument which follows can be applied 
the first time to conclude the existence of a least
positive eigenvalue $\lambda_1$. Iterating this procedure,
a sequence of eigenvalues $\lambda_1 \leq \lambda_2 \leq
\lambda_3 \leq \cdots \leq \lambda_n$ is generated with
corresponding eigenfunctions representative of each
subspace $\omega_1, \, \omega_2, \cdots, \, \omega_n$
for $\Delta$. These can be orthonormalized as well. The
set of eigenfunctions $\{ \omega_1, \cdots, \omega_n \}$
span a subspace of $(H^p)^{\perp}$ which will be called
$V_n$. Clearly both $G$ and $\Delta$ map $(H^p \oplus V_n)^{\perp}
\rightarrow ( H^p \oplus V_n)^{\perp}$. Define $\mu_{n+1}$
to be
$$
\mu_{n+1} = \sup_{\substack{ \left\| \beta \right\| =1\\
\beta \in (H^p \oplus V_n)^{\perp}}} \left\| G \beta \right\|.
\eqno(15)
$$
It will be shown that $\lambda_{n+1} = \mu_{n+1}^{-1}$ is
an eigenvalue of the Laplacian $\Delta$. Moreover,
$$
\mu_{n+1} = \sup_{\substack{\left\| \beta \right\|=1\\
\beta \in (H^p \oplus V_n)^{\perp}}} \left\| G \beta \right\|
\leq \sup_{\substack{\left\| \beta \right\|=1\\ \beta \in
(H^p \oplus V_{n-1})^{\perp}}} \left\| G \beta \right\| = \mu_{n},
\eqno(16)
$$
on account of the fact that $(H^p \oplus V_n)^{\perp}
\subset (H^p \oplus V_{n-1})^{\perp}$. Inequality (16)
will be used to conclude that $\lambda_{n+1} \geq \lambda_n$
at the end. Let $\{ \beta_j \} \in (H^p \oplus V_n)^{\perp}$
be a maximizing sequence for $\mu_{n+1}$, in other words
$\left\| \beta_j \right\| =1$ and $\left\| G \beta_j \right\|
\rightarrow \mu_{n+1}$ as $j \raro \infty$. 
The first claim to be shown that
will be subsequently used is that $\left\| G^2 \beta_j - \mu_{n+1}^2
\beta_j \right\| \rightarrow 0$. Since $G$ is self-adjoint,
it follows that
$$
\left\| G^2 \beta_j - \mu_{n+1}^2 \beta_j \right\|^2
= \left\| G^2 \beta_j \right\|^2 -2 \mu_{n+1}^2
( G^2 \beta_j, \beta_j) + \mu_{n+1}^4
$$
$$
\leq \mu_{n+1}^2 \left\| G \beta_j \right\|^2 - 2 \mu_{n+1}^2
\left\| G \beta_j \right\|^2 + \mu_{n+1}^4 = \mu_{n+1}^4 -
\mu_{n+1}^2 \left\| G \beta_j \right\|^2 \rightarrow 0,
\quad j \rightarrow \infty.
\eqno(17)
$$
Moreover from (17) it follows that  $\left\| G \beta_j
- \mu_{n+1} \beta_j \right\| \rightarrow 0$. Define
$\psi_j = G \beta_j - \mu_{n+1} \beta_j$, then
$$
( \psi_j , G^2 \beta_j - \mu_{n+1}^2 \beta_j )
= ( \psi_j, G (\psi_j + \mu_{n+1} \beta_j)- \mu_{n+1}^2 \beta_j )
= ( \psi_j, G \psi_j + \mu_{n+1} ( G \beta_j - \mu_{n+1} \beta_j))
$$
$$
= ( \psi_j , G \psi_j + \mu_{n+1} \psi_j) =
( \psi_j, G \psi_j) + \mu_{n+1} \left\| \psi_j \right\|^2
\eqno(18)
$$
$$
\geq \mu_{n+1} \left\| \psi_j \right\|^2 >0.
$$
The last inequality follows from the fact that $G$ is a
positive operator, as noted in (14). Since the left-hand side of 
(18) approaches zero as $j \rightarrow \infty$, these
inequalities force $\left\| \psi_j \right\|^2 \rightarrow 0$
as $j \rightarrow \infty$ as well.

Now $\left\| G \beta_j \right\|$ is bounded along with
$\left\| \Delta G \beta_j \right\| = \left\| \beta_j \right\|$.
Therefore, by Proposition 2, there must be a subsequence
of these $\beta_j$, which may as well be called $\{ \beta_j \}$,
such that $\{ G \beta_j \}$ is a Cauchy sequence.
Consequently, a linear functional $L$ can be defined on
$E^p (M)$ by means of the following limit
$$
L ( \sigma ) = \lim_{j \rightarrow \infty} \, \mu_{n+1} \,
( G \beta_j, \sigma ),
\qquad
\sigma \in E^p (M).
\eqno(19)
$$
Then for $u \in ( H^p \oplus V_n)^{\perp}$, next replace
$\sigma$ in (19) by $( \Delta - \frac{1}{\mu_{n+1}}) u$.
Since $\Delta$ is a self-adjoint operator, the definition
in (19) takes the following form
$$
L (( \Delta - \frac{1}{\mu_{n+1}}) u ) = \lim_{j \raro \infty}
\mu_{n+1} ( G \beta_j, \Delta u - \frac{1}{\mu_{n+1}} u)
$$
$$
= \lim_{j \raro \infty} \, \mu_{n+1} ( \Delta G \beta_j, u)
- \mu_{n+1} ( \frac{1}{\mu_{n+1}} G \beta_j, u)
= \lim_{j \raro \infty} \mu_{n+1} ( \frac{1}{\mu_{n+1}}
( \mu_{n+1} \beta_{j} - G \beta_j ), u)
$$
$$
= \lim_{j \raro \infty} (( \mu_{n+1} \beta_j - G \beta_j ), u) =0.
$$
Therefore, this $L$ is a non-trivial, weak solution of the
equation
$$
( \Delta - \frac{1}{\mu_{n+1}} ) u =0.
\eqno(20)
$$
Setting $\lambda_{n+1}  = \mu_{n+1}^{-1}$, Proposition 1
can be applied to this. For $u \in E^p (M)$ and $L$ a 
weak solution of the equation $\Delta u - \lambda_{n+1} u=0$,
there then exists an $\omega_{n+1}  \in E^p (M)$ such that
$L (\beta ) = ( \omega_{n+1}, \beta)$ for every $\beta \in E^p (M)$.
Consequently, this $\omega_{n+1}$ satisfies the equation 
$\Delta \omega_{n+1} = \lambda_{n+1} \omega_{n+1}$. Finally,
this implies that $\lambda_{n+1}$ so defined is an eigenvalue of
$\Delta$.  $\clubsuit$

To finish it might be of interest to formulate a very 
useful result which takes advantage of many of the facts and 
results which have been produced up to now.

{\bf Theorem 7.} Let $\lambda_1 \leq \lambda_2 \leq \cdots$
be the eigenvalues of $\Delta$ on $E^p (M)$ such that each
eigenvalue is included as many times as the dimension of its
eigenspace. Let $\{ \omega_k \}$ be a corresponding sequence
of eigenfunctions. For any $\alpha \in E^p (M)$, 
$$
\lim_{n \raro \infty} \left\| \alpha - \sum_{k=1}^n
 \, ( \alpha, \omega_k) \, \omega_k \right\| =0.
\eqno(21)
$$

{\bf Proof:} Let $m$ be the dimension of the space $H^p$.
There exists $\beta \in (H^p)^{\perp}$ such that
$$
G \beta = \alpha - \sum_{k=1}^{m} \, ( \alpha,
\omega_k ) \, \omega_k,
$$
as in Lemma 1. The set $\{ \omega_i \}$ can be orthonormalized,
which means that $( \omega_i, \omega_j ) =0$ when $i \neq j$ 
and consequently,
$$
\left\| \alpha - \sum_{k=1}^n \, ( \alpha, \omega_k ) \, \omega_k 
\right\| = \left\| \alpha -  \sum_{k=1}^m ( \alpha, \omega_k )
\omega_k - \sum_{k=m+1}^m \, ( \alpha, \omega_k ) \, \omega_k 
\right\|
$$
$$
= \left\| G \beta - \sum_{k=m+1}^n \, ( G \beta + \sum_{j=1}^m
( \alpha, \omega_j ) \omega_j, \omega_k) \, \omega_k \right\|
= \left\| G \beta - \sum_{k=m+1}^n \, ( G \beta, \omega_k) \,
\omega_k \right\|
= \left\| G ( \beta - \sum_{k=m+1}^n \, ( \beta, \omega_k ) \omega_k
\right\|,
$$
for $n>m$. However, for any $v \in ( H^p \oplus V_n )^{\perp}$,
from Theorem 6, it is known that $\left\| Gv \right\| \leq \mu_{n+1}
\left\| v \right\|$. It then follows that
$$
\left\| G ( \beta - \sum_{k=m+1}^n \, ( \beta, \omega_k ) \omega_k )
\right\| \leq \frac{1}{\lambda_{n+1}} \left\| \beta - \sum_{k=m+1}^n \,
( \beta, \omega_k) \omega_k \right\| \leq \frac{1}{\lambda_{n+1}}
\left\| \beta \right\|.
$$
Since $\left\| \beta \right\|$ is a fixed number, by the results
developed pertaining to the eigenvalues, the right-hand
side goes to zero as $n \raro \infty$. $\clubsuit$

\vspace{2cm}
\begin{center}
{\bf References.}
\end{center}

\noindent
{\bf 1.} S. Goldberg, Curvature and Homology, Dover, NY, 1970. \\
{\bf 2.} S. Morita, Geometry of Differential Forms, AMS Translations
of Mathematical Monographs, Vol. 209, Providence, RI, 2001.  \\
{\bf 3.} S. S. Chern, Topics in Differential Geometry, Institute for
Advanced Study Notes, (1951).  \\
{\bf 4.} S. S. Chern, W. H. Chen and K. S. Lam, Lecture Notes on
Differential Geometry, University Mathematics vol. 1, World Scientific,
Singapore, 1999.  \\
{\bf 5.} F. Warner, Foundations of Differentiable Manifolds and Lie Groups,
Springer, NY, 1983.  \\
{\bf 6.} G. de Rham, Vari\'et\'es Differentiables, Hermann, 1955.  \\

\end{document}